\documentstyle[twoside,amssymb,12pt]{article}
\newcommand{\p}{\partial}

\newtheorem{prop}{Proposition}[section]

\newtheorem{theo}[prop]{Theorem}

\newtheorem{coro}[prop]{Corollary}

\title{\sc New families of  conservative
systems on $S^2$ possessing
an integral of fourth degree in momenta} 

\author{{\sc Elena N. Selivanova}\footnote{Supported by DAAD.} \\ }
\begin{document}
\date{}

\maketitle
\thispagestyle{empty}

\section{\bf Introduction}

\noindent
There  is a well-known example of integrable 
conservative system on $S^2$,  the case of Kovalevskaya \cite{Kov} 
in the dynamics of a rigid body, possessing an integral
 of fourth degree in momenta, see \cite{AKN}. 
Using the explicit formular for the  total energy of this 
system
$$H=\frac{du_1^2 + du_2^2 + 2du_3^2}{2u_1^2 + 2u_2^2 + u_3^2} - u_1,$$
where $S^2$ is given by $ u_1^2 + u_2^2 + u_3^2=1$, (see 
 \cite{BKF}), we can rewrite $H$ in polar coordinates as follows
\begin{equation}
H=\gamma_1(r)(r^2d\varphi^2 + dr^2) - \gamma_2(r)\cos \varphi
\label{formh}
\end{equation}
where $\gamma_1$, $\gamma_2$ are some functions 
and $\gamma_2(r)\ne 0$ for $0<r<\infty$.

Goryachev  proposed in \cite{Gor} a  
family of examples of conservative systems on $S^2$
possessing an integral of fourth
degree in momenta. 
In polar coordinates the total energy can be written as follows
\begin{equation}
H=\rho_0(r, B_1, B_2)(r^2d\varphi^2 + dr^2) - \rho_1(r)(
B_1\sin2\varphi + B_2\cos2\varphi) -
\rho_2(r) \cos\varphi
\label{gor}
\end{equation}
where $B_1$, $B_2$ are constants and 
$\rho_0$, $\rho_1$, $\rho_2$ are some  functions. 
It has been shown in  \cite{Gor} that this family reduced to 
the case of Kovalevskaya when $B_1=B_2=0$.

We say that a conservative system on $S^2$ is {\it smooth} 
if the Hamiltonian is a sum of a smooth Riemannian metric on $S^2$ 
(kinetic energy)
and a smooth function $U$ on $S^2$ (potential energy or simply potential).

In this sence above examples of Kovalevskaya and Goryachev
are smooth conservative systems on $S^2$.

In this paper we proposed new examples of {\it smooth}
 conservative systems on $S^2$ 
possessing an integral of {\it fourth} degree in momenta.

In our examples we use  the solution $\psi_0$ of the 
following initial value problem
\begin{equation}
\psi'''\psi' + 2\psi''^2 - 3\psi^2=0, \ \ \psi(0)=0, \psi'(0)=1, \psi''(0)=0. 
\label{inipr}
\end{equation}
We set $\Psi_0(r)=\psi_0(\log r)$, $\Psi_1(r)=\psi'_0(\log r)$, 
$\Psi_2(r)=\psi''_0(\log r)$.
\medbreak
We prove the following theorems.

\begin{theo} The Hamiltonian
\begin{equation}
H=\frac{1}{\Psi_1^2(r)}
\left(d\varphi^2+\frac{dr^2}{r^2}\right) - 
\Psi_1^2(r)(\Psi_2(r)-\Psi_0(r))\cos\varphi
\label{s1}
\end{equation}
defines a smooth conservative system on $S^2$
in polar coordinates. This system  possesses
 an integral of fourth degree and does not possess an integral 
quadratic or linear in momenta.

Hamiltonian (\ref{s1}) cannot be obtained from (\ref{formh}) and  (\ref{gor})
by a change of variables on $S^2$.
\label{main-th1}
\end{theo} 

\begin{theo} There is a constant $p_0$ such that for any 
 $p>p_0$ the Hamiltonians
\begin{equation}
H_p=\frac{\Psi_1^2(r)-\Psi_0^2(r)+p}
{\Psi_1^2(r)}\left(d\varphi^2+\frac{dr^2}{r^2}\right)
-\frac{\Psi_1^2(r)(\Psi_2(r)-\Psi_0(r))\cos\varphi}{\Psi_1^2(r)-\Psi_0^2(r)+p}
\label{s2}
\end{equation}
define  a  one-parameter family of 
 smooth conservative system on $S^2$. These systems 
possess  an integral of fourth degree and does not possess an 
integral quadratic or linear in momenta.

Hamiltonians (\ref{s2}) cannot be obtained from (\ref{formh}) and  (\ref{gor})
by a change of variables on $S^2$.
\label{main-th2}
\end{theo} 
In \cite{Se} we proposed a new {\it one-parameter} family 
of smooth conservative systems on $S^2$ possessing an integral {\it cubic} in momenta.

\medbreak 
The paper is organized as follows. 
In the first chapter we  consider the initial value problem (\ref{inipr}), 
show the existence of the solution $\psi_0$ on ${\bf R}$ and 
investigate  its asymptotic behaviour at infinity.
In the second chapter we obtain a local criterion 
for the integrability of a geodesic flow with an 
integral of  fourth degree and find a class
of metrics with integrable geodesic flows. Finally, in the third chapter 
we prove Theorem \ref{main-th1} and Theorem \ref{main-th2}.

\section{Existence and asymptotic behaviour of the solution}

In this chapter we consider initial value problem (\ref{inipr}), 
prove the existence of the 
solution $\psi_0$ on $\bf R$ and investigate its 
asymptotic behaviour at infinity.

\begin{theo}
The solution $\psi_0(y)$ of initial value problem (\ref{inipr}) 
exists for any $y\in \bf R$
and, moreover, there are $C^{\infty}$ functions $\nu$ and $\mu$ such that
$$\psi_0'(y)=(\exp y)\nu(\exp(-2y))=\exp(- y)\nu(\exp(2y)),$$
$$\psi_0'^2(y)(\psi_0''(y) - \psi_0(y))=\exp (-y)\mu(\exp(-2y))=-
(\exp y)\mu(\exp(2y))$$
and $\nu(t)>0$ for any $t$.
\label{asy-th}
\end{theo}

\noindent
{\em Proof.} Initial value problem (\ref{inipr})
has a unique solution $\psi(y)=\psi_0(y)$
which is positive on $(0,\varepsilon)$ and negative
on $(-\varepsilon, 0)$ for  sufficiently small 
$\varepsilon$.
\medbreak
\noindent
Let us consider the case $y>0$. 

Let $R(y)=\ln \psi(y)$ and $q(y)=R'(y)$.
Then the differential equation from  (\ref{inipr}) is 
equivalent to 
 the following 
system of differential equations of the first order:
\begin{equation}
\dot q=p, \ \
\dot p=\frac{1}{q}\left (-3q^4 - 7q^2p - 2p^2 + 3
\right).
\label{sode}
\end{equation}
System  (\ref{sode}) has two singular points:the 
knots $p=0, q=\pm 1$.

Since  system (\ref{sode}) is symmetric with respect to $q\mapsto -q$, 
$y\mapsto -y$,
 it  suffices to consider the case $q>0$.

Show that  orbit $\Gamma$ of (\ref{sode}), related
to $\psi_0(y)$ for $y\ge 0$, converges to the point 
$q=1, p=0$ as $y\to +\infty$.

For $\psi_0(y)$ we have $R(y)\to -\infty$ as $y\to 0+$ and 
then $q(y)=R'(y)\to +\infty$ as $y\to 0+$.

Write $\Gamma=\{(q,p)| p=-q^{2} + qg(q^{-1})\}$ and 
denote $s=q^{-1}$.

By computation we obtain  from (\ref{sode})
the following differential equation for $g$:
$$g'(1-gs) + 3g^2 - 3s^2=0$$
where the initial condition is the following
$$\lim_{s\to 0+}{g(s)}=\lim_{q\to +\infty}{\frac{p+q^2}{q}}
=\frac{\psi_0''(0)}{\psi_0'(0)}=0.$$
So, we have rewritten (\ref{inipr}) as the 
following initial value problem
\begin{equation}
g'(1-gs) + 3g^2 - 3s^2=0, \ \ g(0)=0
\label{inig}
\end{equation}
The solution of (\ref{inig}) is
$g(s)=s^3$. So,
\begin{equation}
\Gamma=\{(q,p)| p=-q^2 + \frac{1}{q^2}\}
\label{p(q)}
\end{equation}
Since $p(1)=0$ in (\ref{p(q)}), orbit $\Gamma$ converges to the point $q=1, p=0$.
\medbreak
Consider the following linear system of differential equations, related to (\ref{sode}).
\begin{equation}
\dot q=p, \ \
\dot p=-12(q-1) - 7p.
\label{lins}
\end{equation}
The eigenvalues of (\ref{lins}) are equal to $-3$ and $-4$.
Therefore (\ref{sode}) is $C^{\infty}$-conjugate to (\ref{lins})
in a neighborhood of $q=1, p=0$, see \cite{Har}.
Taking into account (\ref{p(q)}), we conclude that there exists
 $C^{\infty}$-functions $\xi_1$, $\xi_2$ such that 
the solution of (\ref{sode}), related to $\psi_0(y)$, $y>0$,
can be written as
$$q=1 + \xi_1(\exp(-4y)),$$
$$p=\xi_2(\exp(-4y))$$
where $\xi_1(0)=\xi_2(0)=0$.

So, we get
$$\psi_0(y)=(\exp y)\beta(\exp(-2y))$$
where 
$$\beta(u)=\exp(-\frac{1}{2}\int_{}^{u}{\frac{\xi_1(t^2)}{t}dt}).$$
Thus, $\beta$ is of class $C^{\infty}$. 

Since $$\psi_0'(y)=\psi_0(y)q=(\exp y)\beta(\exp(-2y))(1 + 
\xi_1(\exp(-4y))=(\exp y)\nu(\exp(-2y)),$$
$\nu$ is of class $C^{\infty}$. Since the values of $q$ for orbit $\Gamma$
belong to $(1, +\infty)$ and $\psi'_0(0)=1$, the function $\psi'_0$ is positive
everywhere. Therefore, $\nu(t)>0$ for any $t\ge 0$.

Then we  compute
$$\psi_0'^2(y)(\psi_0''(y) - \psi_0(y))=\psi_0'^2(y)\psi_0(y)(p+ q^2 -1)=$$
$$\exp(-y)\nu^2(\exp(-2y))\beta(\exp(-2y))
\frac{\xi_2(\exp(-4y)) + 2\xi_1(\exp(-4y)) + \xi_1^2(\exp(-4y))}
{\exp(-4y)}$$ $$=\exp(-y)\mu(\exp(-2y))$$
and, therefore, $\mu$ is of class $C^{\infty}$.

Since $\psi_0(y)=-\psi_0(-y)$, we get 
$$\psi_0'(y)=(\exp y)\nu(\exp(-2y))=\exp(- y)\nu(\exp(2y))$$ and
$$\psi_0'^2(y)(\psi_0''(y) - \psi_0(y))=\exp (-y)\mu(\exp(-2y))=-(\exp y)\mu(\exp(2y)).$$

\hfill $\Box$

\section{\bf A criterion for the integrability}

\noindent 
Consider a metric $ds^2=\Theta(u, v)(du^2 + dv^2)$ in conformal
coordinates $u, v$. It can  also be written as
\begin{equation}
ds^2=\theta(w, \bar w)dwd\bar w
\label{met}
\end{equation}
 where $w=u+ i v$. 
The geodesic flow of $ds^2$ is a Hamiltonian system with Hamiltonian
\begin{equation}
H=\frac{p_wp_{\bar w}}{{4\theta}(w, \bar w)}.
\label{ham}
\end{equation}
A polynomial $F$ in momenta $p_u$, $p_v$ 
can be also written as $$F=\sum_{k=0}^{n}{b_k(w,\bar w)p_w^kp_{\bar w}^{n-k}}$$
where $$b_k=\overline {b_{n-k}}, \ \ k=0,...,n.$$ 

If the polynomial $F$ is an additional integral of the geodesic flow with 
 Hamiltonian (\ref{ham}), then $\{F,H\}=0$ and the following holds
\begin{equation}
\theta \frac{\partial b_{k-1}}{\partial w}+(n-(k-1))b_{k-1}
\frac{\partial \theta}{\partial w}+
\theta \frac{\partial b_{k}}{\partial \bar w}+kb_{k}
\frac{\partial \theta}{\partial \bar w}=0,
\label{systpde}
\end{equation}
where $k=0,... ,n+1 $ and $ b_{-1}=b_{n+1}=0$. 
Substituting   $k=0$ and $k=n+1$ in (\ref{systpde}) we get  immediately 
$$  \frac{\partial}{\partial \bar w}b_0\equiv 0$$ and 
$$\frac{\partial}{\partial  w}b_{n+1}\equiv 0.$$
Note that in another conform coordinate system 
\begin{equation}
\hat{a}_0(\hat{w})=a_0(w)
\left(\frac{d\hat{w}}{dw}\right)^n.
\label{zam}
\end{equation}

\begin{prop}
The geodesic flow of a Riemannian
metric on $S^2$
\begin{equation}
ds^2=\lambda(r^2)(r^2d\varphi^2 + dr^2)
\label{polya}
\end{equation}
does not possess a nontrivial 
integral quadratic in momenta (which does not
depend on $H$ and  linear integrals).
\label{kva}
\end{prop}

\noindent
{\em Proof.} Let us introduce the  conform 
coordinate system $u=r\cos \varphi, v=r\sin \varphi$. Denote $w=u+iv$ and 
consider $\hat{w}=w^{-1}$ for $ds^2$. 

Let $F=a_0(w,\bar w)p_w^2+a_1(w,\bar w)p_wp_{\bar w}+a_2(w,\bar w)p_{\bar w}^2$
be an integral of the geodesic of (\ref{polya}). As above, $a_0=a_0(w)$, $a_2=a_2(\bar w)=
\bar a_0(w)$ and $\mbox {Im}\ a_1(w,\bar w)=0$.
Since (\ref{polya}) is a Riemannian
metric on $S^2$, from (\ref{zam}) it follows that 
$a_0(w)$ is a polynomial of degree $n$ where $n\le 4$. 
If $F$ does not depend on
the Hamiltonian and linear integrals, then $a_0(w)$ does not equal to
zero identically.

Then from (\ref{systpde}) when $k=1$ we obtain
$$\lambda\frac{\p a_0}{\p w} + 2a_0\frac{\p \lambda}{\p w} +
\lambda\frac{\p a_1}{\p \bar w} + a_1\frac{\p \lambda}{\p \bar w}=0.$$
This is equivalent to the following condition
$$\mbox {Im}\ \frac{\p}{\p w}
\left(\lambda\frac{\p a_0}{\p w} + 2a_0\frac{\p \lambda}{\p w} \right)=0$$ 
where $a_0(w)=A_0w^4 + A_1w^3 + A_2w^2 + A_3w + A_4$ and  $A_0$,
 $A_1$,  $A_2$,  $A_3$,  $A_4$ are some constants such that at least 
one of them does not equal to zero. 

Denote $w\bar w$ as $t$. So, we get the following conditions for $\lambda(t)$
$$\mbox{Im}\ (w^2(12A_0\lambda +12A_0\lambda't+2A_0\lambda''t^2)+
\bar w^2(2A_4\lambda''))+$$
$$\mbox{Im}\ (w(6A_1\lambda +9A_1\lambda't+2A_1\lambda''t^2)+
\bar w(3A_3\lambda' + 2A_3\lambda''t))+$$
$$\mbox{Im}\ (A_2(2\lambda + 6\lambda't + 2\lambda''t^2))\equiv 0.$$

So, we have to consider the following linear differential equations
\begin{equation}
\epsilon_1(6\lambda+6\lambda't+\lambda''t^2)=\epsilon_2(\lambda''),
\label{dif1}
\end{equation}
\begin{equation}
\epsilon_1(6\lambda+9\lambda't+2\lambda''t^2)=\epsilon_2(3\lambda'+
\lambda''t),
\label{dif2}
\end{equation}
\begin{equation}
\lambda + 3\lambda't + 2\lambda''t^2=0
\label{dif3}
\end{equation}
where $\epsilon_1$, $\epsilon_2$ are constants equal to $-1$, $0$ or $1$.
In another case consider the rescaling 
$t\mapsto Dt$, $D$ is a real constant,  of a solution
(\ref{dif1}) or (\ref{dif2})). 

It is easy to see that for any smooth solution of (\ref{dif3}) $\lambda(0)=0$
and then $ds^2$  is not a Riemannian metric on $S^2$. (We get that 
$\mbox{Im}\ A_2=0$.)

It is also easy to see that any solution of the equations  (\ref{dif1}), 
 (\ref{dif2}) if $\epsilon_1=0$ or $\epsilon_2=0$  does not define 
by (\ref{polya}) a 
Riemannian metric on $S^2$.

Solutions of (\ref{dif1}) if $\epsilon_1=1$ and $\epsilon_2=1$
have the following form
$$\lambda(t)=C_1\frac{1+t^2}{(-1+t^2)^2} 
+C_2\frac{t}{(-1+t^2)^2}, C_1, C_2 - const.$$
Solutions of (\ref{dif1}) if $\epsilon_1=1$ and $\epsilon_2=-1$
have the following form
$$\lambda(t)=C_1\frac{1-t^2}{(1+t^2)^2} 
+C_2\frac{t}{(1+t^2)^2}, C_1, C_2 - const.$$
We got it noting that solutions of 
(\ref{dif1}) if $\epsilon_1=1$ and $\epsilon_2=1$ can
be obtained from the solutions of (\ref{dif1}) if 
$\epsilon_1=1$ and $\epsilon_2=-1$ by complex transformation 
$t\to it$.

It is easy to note that if $\lambda(0)=1$ there is only one smooth
solution of 
 (\ref{dif2}) if $\epsilon_1=1$ or $\epsilon_2=\pm 1$. These solutions 
correspond to the metrics of curvature $-1$ if $\epsilon_2=1$ 
and to the metrics of curvature $-1$ if $\epsilon_2=-1$. So, 
if $\epsilon_1=1$ and
$\epsilon_2=1$
any solution of  (\ref{dif2}) which is smooth in zero   has the following form
$$\lambda(t)=\frac{C_1}{(-1 +t)^2}, C_1 - const$$
and  if $\epsilon_1=1$ and
$\epsilon_2=-1$
any solution of  (\ref{dif2}) which is smooth in zero has the following form
$$\lambda(t)=\frac{C_1}{(1 +t)^2}, C_1 - const.$$
Thus, we may see that either $\mbox {Im}\  A_2=0$ and $A_i=0, i\ne 2$
or $$\lambda(t)=\frac{C_1}{(-1 +Dt)^2}, C_1, D - const$$ and, therefore,
\begin{equation}
ds^2=\frac{C_1}{(1+Dr^2)^2}\left(r^2d\varphi^2 +dr^2\right), 
\label{ccur}
\end{equation}
where $C_1, D - const$, i.e. 
$ds^2$ is a metric of constant positive curvature.

Remember that the geodesic flow of a metric (\ref{polya}) possesses the 
linear integral $F_1=iwp_w-i\bar wp_{\bar w}$.
If $\mbox {Im}\  A_2=0$ and $A_i=0, i\ne 2$, we will consider 
the integral $F_2=F+A_2F_1$ which is equal, obviously, $\hat{C}H$ where
 $\hat{C}$ is a constant. Therefore, $F$ depends on $H$ and $F_1$.

So, if the geodesic flow of a metric (\ref{polya}) on $S^2$
possesses a nontrivial quadratic integral, it can be only a metric
of constant positive curvature. But it is well known that 
the geodesic flows of metrics of constant curvature possess 
two independent linear integrals. (Indeed, the geodesic flows
with the Hamiltonians $H=l'^2(y)(p_x^2 + p_y^2)$ where $l''(y)=l(y)$
possess the integrals $F_1=p_x$ and $\hat{F_1}=l(y)\cos xp_x - l'(y)\sin xp_y$.)
So, in this case any quadratic integral is also trivial, i.e. depends on 
the Hamiltonian and linear integrals.
\hfill  $\Box$ 

\begin{prop} 
The geodesic flow of a metric (\ref{polya}) on $S^2$
possesses two independent linear  integrals if and only if it 
 has form (\ref{ccur}), i.e. it is 
a metric of constant positive curvature. 
\label{kva2}
\end{prop}

\noindent
{\em Proof.} As mentioned above the geodesic flow of 
a metric (\ref{polya}) on $S^2$
possesses the linear integral $F_1=p_{\varphi}$.

In the coordinates $w, \bar w, p_w, p_{\bar w}$ where 
$w=r\cos \varphi + ir\sin \varphi$ it can be written as 
follows $F_1=iwp_w - i\bar wp_{\bar w}$.
(Indeed, in the coordinates $z=\varphi + iy$, $y=\log r$ we 
have $F_1=p_z+p_{\bar z}$ and, clearly, $w=\exp{iz}$.)

Assume that a metric $ds^2$ of form (\ref{polya}) on $S^2$
possesses another linear integral $F_2=b_0(w)p_w + \bar b_0(w)p_{\bar w}$ 
which is  independent of $F_1$ and, therefore, the function 
$b_0(w)$ does not equal to $iwC_3$ where $C_3$ is a real constant.

Let us consider then the function $F_2^2=b_0^2(w)P_w + 2|b_0|^2p_wp_{\bar w}
+ \bar b_0(w)p_{\bar w}^2$ which is also an integral of the geodesic flow 
of $ds^2$.

Compare now $b_0^2(w)$ with $a_0(w)$ from Proposition \ref{kva}.
Since $b_0(w)\ne iwC_3$ in our case, then at least one coefficient 
$A_i$, $i\ne 2$ (see Proposition \ref{kva}) does not equal to zero.
In the same way as in Proposition \ref{kva} we may show that 
$ds^2$ has form (\ref{ccur}),
i.e. 
$ds^2$ is a metric of constant positive curvature.

\hfill  $\Box$

Remember that conform coordinates $x,y$ of a metric $ds^2$ are 
Liouville if $ds^2=(f(x)+h(y))(dx^2 +dy^2)$ for some functions $f$, $h$.
Due to Darboux \cite{Dar} Liouville coordinates exist if and only if
the geodesic flow of $ds^2$ possesses an integral quadratic in momenta
and moreover if one of the function $f$ or $h$ is constant then
the geodesic flow possesses an integral linear in momenta.
So, we may formulate the following 
\begin{coro}
Liouville  coordinates $\varphi, y=\log r$, related
to polar coordinates $\varphi, r$
of a metric (\ref{polya}) 
on $S^2$ are unique up to shifts and the transform $y\to -y$.
\label{typ}
\end{coro}

\noindent
{\em Proof.} Assume that a metric $ds^2$ 
on $S^2$ can be written in two different forms (\ref{polya}) 
$$ds^2=\lambda_1(r^2)(r^2d\varphi^2 + dr^2)$$
and 
$$ds^2=\lambda_2(\tilde r^2)(\tilde r^2d\tilde\varphi^2 + d\tilde r^2)$$
where $\tilde r\ne Dr^{\pm 1}$, $D - const$. This means that 
the geodesic flow of $ds^2$ possesses two independent linear integrals.
So, from Proposition \ref{kva2} it follows that 
$\tilde r\ne Dr^{\pm 1}$, $D - const$ and, therefore,
Liouville coordinates $\varphi, y=\log r$ 
 related
to polar coordinates $\varphi, r$
of $ds^2$ on $S^2$  are unique up to shifts and the transform $y\to -y$.

\hfill  $\Box$ 

We must note that Proposition \ref{kva} follows 
from the results of Kolokol'tsov
published in his Ph.D. thesis (in Russian) but we give here a complete 
proof.

One may show that if the polynomial in momenta integral $F$  
of the geodesic flow of 
(\ref{met})  is independent of the Hamiltonian 
and an integral of smaller degree, then there is a conformal coordinate system $z=z(w)$ 
of this  metric such that the coefficients of $F$ for $p_z$ and $p_{\bar z}$
 equal $1$ identically.

\begin{theo} Let $ds^2=\lambda(z, \bar z)dzd\bar z$ be a metric
such that there exists a function $f:{\bf R}^2\mapsto \bf R$, satisfying the 
following conditions 
\begin{equation}
\lambda(z,\bar z)=\frac{\partial^2f}{\partial  z\partial \bar z}, \ \ 
\mbox{Im}\ \left(\frac{\partial^4f}{\partial  z^4}
\frac{\partial^2f}{\partial  z\partial \bar z} + 
3\frac{\partial^3f}{\partial  z^3}\frac{\partial^3f}
{\partial  z^2\partial \bar z} + 2\frac{\partial^2f}{\partial  z^2}
\frac{\partial^4f}{\partial  z^3\partial \bar z}\right)=0.
\label{eqpde}
\end{equation}
Then the geodesic flow of $ds^2$ possesses an integral of fourth degree in momenta.
\medbreak
If the geodesic flow of a metric $ds^2$ possesses an integral 
which is a polynomial 
 of  fourth degree in momenta and it does not depend on 
the Hamiltonian and an  integral of smaller degree then there exist
  conformal coordinates $x, y$ and a function   $f:{\bf R}^2\mapsto \bf R$ such that 
$ds^2=\lambda(z, \bar z)dzd\bar z$ where $z=x+i y$ and (\ref{eqpde}) holds.
\label{loc-th}
\end{theo} 

\noindent
{\em Proof.} Suppose that (\ref{eqpde}) holds. We will now construct an 
integral $$F=\sum_{k=0}^{n}a_k(z, \bar z)p_z^kp_{\bar z}^{n-k}$$ 
of the corresponding geodesic flow. 

Put ${a}_0=1$ and ${a}_n=1$. Equation (\ref{systpde}) then has  the following form
\begin{equation}
4\frac{\partial \lambda}{\partial  z}=-\frac{\partial (a_1\lambda)}{\partial \bar z},
\label{ur1}
\end{equation}
\begin{equation}
\frac{\partial (a_1\lambda^3)}{\partial  z}=-\lambda
\frac{\partial (a_2\lambda^2)}{\partial \bar z},
\label{ur2}
\end{equation}
\begin{equation}
a_2=\overline{a_2}.
\label{ur3}
\end{equation}
Consider $$a_1=-4\frac{\partial^2f}{\partial  z^2}
\left(\frac{\partial^2f}{\partial  z\partial \bar z}\right)^{-1}.$$ Then 
(\ref{ur1}) holds.

Show that for $$g=\frac{\partial (a_1\lambda^3)}{\partial  z}\lambda^{-1}$$ it holds
 $$\mbox {Im}\ 
\frac{\partial g}{\partial  z}=0.$$ 
Indeed, by  computation we obtain
$$\mbox {Im}\ 
\frac{\partial g}{\partial  z}=-4\
\mbox{Im} \left(\frac{\partial^4f}{\partial  z^4}
\frac{\partial^2f}{\partial  z\partial \bar z} + 
3\frac{\partial^3f}{\partial  z^3}\frac{\partial^3f}
{\partial  z^2\partial \bar z} + 2\frac{\partial^2f}{\partial  z^2}
\frac{\partial^4f}{\partial  z^3\partial \bar z}\right)=0.$$
So, it follows  that there is a function $h:{\bf R}^2\mapsto \bf R$ such that 
$$\frac{\partial h}{\partial \bar z}=g.$$ 

Consider $a_2=-\frac{h}{\lambda^2}$. Then $a_2: {\bf R}^2\mapsto \bf R$  and, moreover,
$${\frac{\partial (a_1\lambda^3)}{\partial  z}}{\lambda}^{-1}=-
\frac{\partial (a_2\lambda^2)}{\partial \bar z},$$
i.e. conditions (\ref{ur2}) and (\ref{ur3}) hold also. 

Put then $a_3=\bar a_1$. So, we have an  integral of fourth degree
in momenta of the geodesic flow of $ds^2$. 
\bigbreak
Prove now the second statement of the theorem.
As mentioned above, in this case there exists a conformal 
coordinate system of 
$ds^2$ where $a_0=a_4\equiv 1$. Note that this system is unique up to a shift.
So, in this coordinate system (\ref{ur1})-(\ref{ur3}) hold. From (\ref{ur1}) it follows
that there exists a function $V$ where 
$$\frac{\partial V}{\partial  \bar z}=\lambda$$ and 
$$a_1\lambda=-4\frac{\partial V}{\partial   z}.$$
From (\ref{ur2}) and (\ref{ur3}) it follows immediately that
$$\mbox{Im}\ \frac{\partial }{\partial z}\left(\frac{\partial (a_1\lambda^3)}
{\partial z}\lambda^{-1} \right)=0.$$
On the other hand
$$\mbox{Im}\ \frac{\partial }{\partial z}
\left(\frac{\partial (a_1\lambda^3)}{\partial z}\lambda^{-1} \right)=
-4\ \mbox{Im}\ \frac{\partial }{\partial z}
\left( \frac{\partial^2 V }{\partial z^2}
\frac{\partial V }{\partial \bar z} + 2\frac{\partial^2 V }{\partial z\partial \bar z}
\frac{\partial V }{\partial  z}\right)=0.$$
Taking into account that $\mbox{Im}\ \frac{\partial V}{\partial \bar z}=\mbox{Im}\ \lambda=0$, we
conclude that there exists a function $f:{\bf R}^2\mapsto \bf R$ such that
$$\frac{\partial f}{\partial z}=V.$$
Thus, $$\lambda=\frac{\partial^2 f}{\partial z\partial \bar z}$$ and
$$0=\mbox{Im}\ \frac{\partial }{\partial z}\left( \frac{\partial^3 f}{\partial z^3}
\frac{\partial^2 f}{\partial z\partial \bar z} +
2\frac{\partial^3 f}{\partial z^2\partial \bar z}
\frac{\partial^2 f}{\partial z^2}\right)=$$
$$=\mbox{Im}\ \left(\frac{\partial^4f}{\partial  z^4}
\frac{\partial^2f}{\partial  z\partial \bar z} + 
3\frac{\partial^3f}{\partial  z^3}\frac{\partial^3f}
{\partial  z^2\partial \bar z} + 2\frac{\partial^2f}{\partial  z^2}
\frac{\partial^4f}{\partial  z^3\partial \bar z}\right).$$

\hfill  $\Box$ 

Equation (\ref{eqpde}) in some other form has been written in 
\cite{Hall} but in order to show when the corresponding integral is nontrivial
we needed to give here a complete proof.

Using this criterion we may find families of metrics with
integrable geodesic flows.

We will look for a class of  solutions of (\ref{eqpde})
and show that in this class of  solutions the problem of 
integrability can be reduced to an ordinary differential equation.  

Denote $$\frac{\partial^2
 f}{\partial z^2}=\left (A(x,y) + i B(x,y) \right).$$
By  computation we obtain
$$\mbox {Im}\ \frac{\partial^4 f}{\partial z^4}=\frac{1}{4}
\left(B_{xx} - B_{yy} -2A_{xy}\right),$$
$$\mbox {Im}\ \frac{\partial^3f}{\partial  z^3}\frac{\partial^3f}
{\partial  z^2\partial \bar z}=\frac{1}{2}\left(A_xB_x + A_yB_y\right),$$
$$\mbox {Im}\ \frac{\partial^2f}{\partial  z^2}
\frac{\partial^4f}{\partial  z^3\partial \bar z}=
\frac{1}{4}\left( B(A_{xx} + A_{yy}) + A(B_{xx} + B_{yy}) \right).$$
Thus,  (\ref{eqpde}) can be written  as follows
\begin{equation}
6(A_yB_y + A_xB_x) + 2(B(A_{xx}+A_{yy}) + A(B_{xx}+B_{yy}))
+(B_{xx}-B_{yy}-2A_{xy})\lambda =0.
\label{AB}
\end{equation}
Consider now the solutions of (\ref{eqpde}) of the following form:
\begin{equation}
f(x,y)=\psi(y)\cos x + \xi(y) + d(x^2 - y^2)
\label{f}
\end{equation}
where $\psi$, $\xi$ are some smooth functions and $d$ is a constant.

Then we obtain
$$
\lambda=\frac{1}{4}\left((\psi''(y)- \psi(y))\cos x + \xi''(y) \right)
$$
and
\begin{equation}
A=-\frac{1}{4}\left( (\psi''(y)+ \psi(y))\cos x + \xi''(y) -4d \right),
 B=\frac{1}{2}\left(\psi'(y)\sin x \right).
\label{A,B}
\end{equation}
Substituting (\ref{A,B}) in (\ref{AB}) we get the 
following condition of integrability
$$6(((\psi''' + \psi')\cos x + \xi''')2\psi''\sin x -
 (\psi'' + \psi)2\psi'\sin x\cos x + $$
$$2(2\psi'\sin x((\psi^{(4)} - \psi)\cos x + \xi^{(4)}) + 
((\psi''+ \psi)\cos x + \xi'' - 4d)2(\psi'''-\psi')\sin x))+$$
$$4(\psi''' + \psi')\sin x((\psi''-\psi)\cos x + \xi'')\equiv 0.$$
Therefore,
\begin{equation}
3\psi''\xi''' + \psi'\xi^{(4)} + 2\psi'''\xi'' = 4d(\psi''' - \psi'),
\label{eqxi}
\end{equation}
\begin{equation}
5\psi'''\psi'' - 6\psi'\psi + \psi'\psi^{(4)}=0.
\label{eqpsi}
\end{equation}
Equation (\ref{eqxi}) can be written as
$$(\psi'\xi''')' + 2(\psi''\xi'')' = 4d(\psi'' - \psi)'.$$
By integrating we get
$$\psi'\xi''' + 2\psi''\xi''=4d(\psi'' - \psi) + d_1$$
where $d_1$ is a constant.
So,  if $d=0$ 
\begin{equation}
\xi''=\frac{d_1\psi(y) + c}{(\psi'(y))^2}.
\label{xi1}
\end{equation} or if $d\ne 0$
\begin{equation}
\xi''=2d\frac{\psi'^2(y) - \psi^2(y) +
 d_1(2d)^{-1}\psi(y) + p}{\psi'(y))^2}.
\label{xi2}
\end{equation}
where $p$, $c$ are some constants.

On the other hand there is an integral of (\ref{eqpsi}):
\begin{equation}
2\psi''^2 - 3\psi^2 + \psi'\psi''' = const.
\label{0}
\end{equation}
So, we may now formulate the following proposition. 
\begin{prop}
For any solution  $\psi$ of  (\ref{0})
  there are
two families of metrics $ds_1^2=\Lambda(x,y)(dx^2 + dy^2)$  where 
\begin{equation}
\Lambda(x,y)=(\psi''(y)- \psi(y))\cos x + 
\frac{d_1\psi(y) + c}{\psi'^2(y)}
\label{fam1}
\end{equation}
or
\begin{equation}
\Lambda(x,y)=
(\psi''(y)- \psi(y))\cos x + 
c\frac{\psi'^2(y) - \psi^2(y) +
 d_1\psi(y) + p}{\psi'^2(y)} 
\label{fam2}
\end{equation}
 where $c$, $d$, $p$ are constants and the corresponding
 geodesic flows possess an integral
of fourth degree in momenta. If $\frac{\p^2\Lambda}{\p x\p y}$ and  
$\frac{\p^2\Lambda}{\p x^2}-\frac{\p^2\Lambda}{\p y^2}$ are not equal 
zero identically, then this integral is nontrivial, i.e. it does not depend
on the Hamiltonian and an integral of  smaller degree.
\label{class-th}
\end{prop}

\noindent
{\em Proof.} As we have shown, for every metric in one of these families
 there exists a function
$f$ of the form (\ref{f}) satisfying (\ref{eqpde}). In fact, 
compare (\ref{fam1}) with (\ref{xi1}) and (\ref{fam2}) with (\ref{xi2})
(in the last case we denote  $4d$ in (\ref{xi2}) as $c$). 

Since every metric in one of these families satisfies the conditions 
of theorem \ref{loc-th}, the integrability follows. 

So, we have only to prove  that the integrals of the geodesic flows of these
geodesic flows are nontrivial. 

According to our construction the integral $F$ of the geodesic flow  
with  Hamiltonian $$H=\frac{p_zp_{\bar z}}{4\Lambda}$$ 
where $\Lambda$ is given by (\ref{fam1}), (\ref{fam2}) has the following form
$$F=p_z^4 + D(x, y)p_z^3p_{\bar z} + E(x, y)p_z^2p_{\bar z}^2 + \bar D(x, y)p_zp_{\bar z}^3 + p_{\bar z}^4$$
where $D(x, y): \bf R^2 \mapsto \bf C$,  $E(x, y): \bf R^2 \mapsto \bf R$.

Thus, if $F$ depends on $H$ and another  integral polynomial of smaller degree, then there is an 
integral $\tilde F$  quadratic in momenta which has the following form
$$\tilde F=p_z^2 + \tilde E_1(x, y)p_z^2p_{\bar z}^2  + p_{\bar z}^2$$
or
$$\tilde F=ip_z^2 + \tilde E_2(x, y)p_z^2p_{\bar z}^2  -ip_{\bar z}^2$$
where $\tilde E_1(x, y), \tilde E_2(x, y) : \bf R^2 \mapsto \bf R$.

Therefore, due to Darboux \cite{Dar} the following holds
$$\frac{\p^2 \Lambda}{\p x\p y}\equiv 0$$ or
$$\frac{\p^2\Lambda}{\p x^2}-\frac{\p^2\Lambda}{\p y^2}\equiv 0.$$

\hfill$\Box$

\section {Smooth metrics on $S^2$}

In this chapter we will prove Theorem \ref{main-th1} and Theorem \ref{main-th2}.

\medbreak
\noindent
{\em Proof of  Theorem \ref{main-th1}.} 
Prove that Hamiltonian system with Hamiltonian (\ref{s1}) is
a smooth system on $S^2$. Using theorem \ref{asy-th}, we may write
$$H=\frac{1}{\nu^2(r^2)}(r^2d\varphi^2 + dr^2) + \mu(r^2)r\cos \varphi=$$
$$\frac{1}{\nu ^2(\tilde r^2)}
(\tilde r^2d\tilde \varphi^2 + d\tilde r^2) - 
\mu(\tilde r^2)\tilde r\cos \tilde \varphi $$
where $\tilde r=\frac{1}{r}$, $\tilde \varphi=-\varphi$.
Since $\nu$, $\mu$ are of class $C^{\infty}$ and $\nu\ne 0$, 
see Theorem \ref{asy-th}, then this system is a smooth system on $S^2$.

Prove that  (\ref{s1}) possesses an integral of fourth degree in momenta.
Consider  metrics (\ref{fam1}) 
 from Proposition \ref{class-th} and 
put in (\ref{fam1}) $\psi(y)=\psi_0(y)$ and  $d_1=0$. 
Then we obtain the following
one-parameter family of metrics:
$$\left((\psi_0''(y) - \psi_0(y))\cos x + c\frac{1}{\psi_0'^2(y)}\right)(dx^2 + dy^2)$$
or in polar coordinates $\varphi=x$, $r=\exp y$:
$$\left((\Psi_2(r) - \Psi_0(r))\cos \varphi + 
c\frac{1}{\Psi_1^2(r)}\right)(d\varphi^2 + \frac{1}{r^2}dr^2)$$
where $c$ is an arbitrary constant. From Proposition \ref{class-th}
we know that any metric in this family  possesses an integral of fourth degree in momenta.
Using the well-known Maupertuis's principle, see \cite{Arn}, we conclude
that (\ref{s1}) possesses an integral of fourth degree in momenta, see also 
 \cite{BaNe}, \cite{BKF}.
\medbreak
So, we have to prove only that this integral nontrivial, i.e. depends on
the Hamiltonian and an integral of smaller degree.
Denote the kinetik energy in (\ref{s1}) as $K=\lambda(y)(d\varphi^2+dy^2)$ and 
the potential as $V$. We will denote the geodesic flow of the metric of $K$ as $K$ also.

Let us assume that  a system of this family has an integral
 which is independent
of  the Hamiltonian (total energy (\ref{s1}))
and which is a polynomial of second degree in momenta (clearly, this assumption includes
 the case of linear integrals). So, there is an  integral $\tilde F$ 
of (\ref{s1}) which is quadratic in momenta.
Thus, $\tilde F= A(p_{\varphi}, p_{y}, \varphi, y) + B(\varphi, y)$ 
where $ A(p_{\varphi}, p_{y}, \varphi, y)$
is a polynomial of  second degree in momenta.
We may write $\{\tilde F, H\}=\{
A(p_{\varphi}, p_{y}, \varphi, y) + B(\varphi, y), 
K +V\}\equiv 0$
and, therefore, $\{A(p_{\varphi}, p_{y}, \varphi, y), K\}\equiv 0$. 

Thus, the geodesic flow of $K$ has an integral which is
a polynomial of  second degree with respect to momenta.

Since the function $\psi$ such that $\psi'(y)=c_1\exp(-y)(1+\exp(2y)), c_1 - const$
does not satisfy (\ref{inipr}) and, therefore, $\Psi_1(r)\ne r^{-1}(1+Dr^2)$, $D$ is
a positive constant, then $K$ is not a metric of constant positive curvature.
So, the quadratic integral $A$ of the geodesic flow of $K$ depends on the Hamiltonian 
$K$ and the integral $p_{\varphi}$. W.l.o.g. we may put $A=p_{\varphi}^2$.

So, we may write
$$\{A,V\}+\{B,K\}=\{p_{\varphi}^2,V\}+\{B,K\}\equiv 0.$$
By computation we obtain
$$\frac{\p B}{\p y}\equiv 0,$$ and
$$\frac{\p B}{\p \varphi}=\frac{\p V}{\p \varphi}\lambda(y).$$
So, we get
$$V=\frac{B(\varphi)+\alpha(y)}{\lambda(y)}$$
for a smooth function $\alpha(y)$. 

In our case $$\lambda(y)=\frac{1}{\psi_0'^{2}(y)}$$ and
$$V(x,y)=(\psi_0''(y)-\psi_0(y))\psi_0'^{2}(y)\cos \varphi$$
and we get then
$\psi_0''(y)-\psi_0(y)\equiv const$, that is not true.
So, there is no quadratic integral of the system given by (\ref{s1}).

\medbreak
In oder to prove that the system given by (\ref{s1}) is really a new example
of the integrability by a polynomial of fourth degree in momenta we have to
compare this case with the known cases of Kovalevskaya and Goryachev.
 First of all, we note that in the family (\ref{gor})
we have to consider only the case $B_1=B_2=0$. As  mentioned above,
for $B_1=B_2=0$ the system given by (\ref{gor}) is simply the case of Kovalevskaya.
So, we may consider only the case of Kovalevskaya.  Then 
comparing the potentials in (\ref{formh}) with the potentials 
in (\ref{s1}) and (\ref{s2}) we see that they are cannot be the same, because
$\gamma_2(r)>0$, $r\in (0, +\infty)$ in (\ref{formh}) but $\Psi_2(1)-\Psi_0(1)=0$.
 
So, this example is really new, i.e. it  cannot be obtained from the known 
integrable coservative systems on $S^2$ corresponding to the cases of Kovalevskaya 
and Goryachev.

\hfill $\Box$

\medbreak
\noindent
{\em Proof of  Theorem \ref{main-th2}.} 
Consider  Hamiltonian (\ref{s2}) and prove that there is a constant $p_0$ that
for any $p>p_0$ it is 
a smooth system on $S^2$. Using theorem \ref{asy-th}, we may rewrite $H_p$
in (\ref{s2}) as 
$$H_p=\frac{\Phi(r^2)+p+1}{\nu^2(r^2)}(r^2d\varphi^2 + dr^2)
 + \frac{\mu(r^2)}{\Phi(r^2)+p+1}r\cos \varphi=$$
$$\frac{-\Phi(\tilde r^2)+p+1}{\nu ^2(\tilde r^2)}
(\tilde r^2d\tilde \varphi^2 + d\tilde r^2) - 
\frac{\mu(\tilde r^2)}{-\Phi(\tilde r^2)+p+1}\tilde r\cos \tilde \varphi $$
where $\tilde r=\frac{1}{r}$, $\tilde \varphi=-\varphi$ and 
$$\Phi(t)=\int_{t}^{1}{\mu(s)\nu^{-1}(s)ds}.$$

Since $\nu$, $\mu$ are of class $C^{\infty}$ and $\nu\ne 0$, 
see Theorem \ref{asy-th},  there are constants 
$$M_1=\max_{[0,1]}\Phi(t), \  M_2=-\min_{[0,1]}\Phi(t)$$
and, therefore, for $p>p_0=\max \{M_1, M_2\} -1$ system with Hamiltonian
(\ref{s2}) is a smooth conservative system on $S^2$.

Prove that the system with energy (\ref{s2}) 
possesses an integral of fourth degree in momenta.
Consider  metrics (\ref{fam2}) 
 from Proposition \ref{class-th} and 
put in (\ref{fam2}) $\psi(y)=\psi_0(y)$ and  $d_1=0$. 
We get the following two-parameter
family of metrics:
$$\left((\psi_0''(y) - \psi_0(y))\cos x + c\frac{\psi_0'^2(y) - \psi_0^2(y) + p}{\psi_0'^2(y)}\right)(dx^2 + dy^2)$$
which can be written also  in polar coordinates $\varphi=x$, $r=\exp y$ as
$$\left((\Psi_2(r) - \Psi_0(r))\cos \varphi + 
c\frac{\Psi_1^2(r) - \Psi_0^2(r) + p}{\Psi_1^2(r)}\right)(d\varphi^2 + 
\frac{1}{r^2}dr^2)$$
where $c$ is an arbitrary constant. It was shown in Proposition \ref{class-th}
that any metric in this family  possesses an integral of fourth degree in momenta.
From the well-known Maupertuis's principle, see \cite{Arn}, it follows
that the system with Hamiltonian (\ref{s2}) possesses an integral of fourth degree in momenta.

\medbreak
In order to prove  that for any $p>p_0$ the system with Hamiltonian (\ref{s2})
does not possess an integral quadratic or linear in momenta and to compare 
with the above cases of Kovalevskaya and Goryachev we can apply the same arguments 
as for the system with Hamiltonian (\ref{s1}), see the proof of Theorem \ref{main-th1}.

\medbreak
We must prove that for different values of parameter $p$ the systems given by (\ref{s2})
cannot be obtained one from another by a change of variables on $S^2$, i.e. by
a diffeomorphismus $\phi :\ S^2\to S^2$.

Let us introduce the coordinates $\varphi, y=\log r$.
From Corollary \ref{typ} it follows that we need to consider only  
the transformations $y\mapsto \pm y + \kappa$ where $\kappa$ is a constant.

Let us rewrite $H_p$ in the coordinates $\varphi, y$:
$$H_p=\frac{\psi_0'^2(y)-\psi_0^2(y) +p}{\psi_0'^2(y)}\left(d\varphi^2 + dy^2\right)-$$
\begin{equation}
(\psi_0''(y)-\psi_0(y))\left(\frac{\psi_0'^2(y)-\psi_0^2(y) +p}{\psi_0'^2(y)}\right)^{-1}
\cos\varphi
\label{Hy}
\end{equation}
If the systems given by (\ref{Hy}) for $p=p_1$ and $p=p_2$ can be obtained one from another
by a transform $y\mapsto \pm y + \kappa$, $\kappa - const$, then there are 
some constants $K_1\ne 0$, $K_2\ne 0$ such that the following holds
$$\frac{\psi_0'^2(y)-\psi_0^2(y) +p}{\psi_0'^2(y)}=$$
\begin{equation}
K_1\frac{\psi_0'^2(\pm y + \kappa )-\psi_0^2(\pm y + \kappa) +p}{\psi_0'^2(\pm y + \kappa)}
\label{K1}
\end{equation}
and 
$$(\psi_0''(y)-\psi_0(y))\left(\frac{\psi_0'^2(y)-\psi_0^2(y) +p}{\psi_0'^2(y)}\right)^{-1}=$$
\begin{equation}
K_2(\psi_0''(\pm y + \kappa)-\psi_0(\pm y  + \kappa))\left(\frac{\psi_0'^2(\pm y + \kappa)-
\psi_0^2(\pm y + \kappa) +p}{\psi_0'^2(\pm y + \kappa)}\right)^{-1}.
\label{K2}
\end{equation}
Thus, we get
\begin{equation}
\psi_0''(y)-\psi_0(y)=K_1K_2(\psi_0''(\pm y + \kappa)-\psi_0(\pm y + \kappa)).
\label{K1K2}
\end{equation}
Put $y=0$ in (\ref{K1K2}) and we get $\psi_0''(\kappa)-\psi_0(\kappa)=0$.
W.l.o.g. we assume that $\kappa \ge 0$.
From the proof of Theorem \ref{asy-th} we know that the function $\psi_0''(y)-\psi_0(y)$, $y>0$
can be written as follows
$$\psi_0''(y)-\psi_0(y)=(\exp{R(y)})(p(q) + q^2 - q), y>0$$
where $p(q)=-q^2+\frac{1}{q^2}$ and $q>1$.

So, we obtain $\psi_0''(y)-\psi_0(y)=\exp{R(y)}(q^{-2}-q)$ and, therefore,
if $\psi_0''(\kappa)-\psi_0(\kappa)=0$ and $\kappa\ne 0$, then $q(\kappa)=1$ that
is impossiblle. Thus, $\kappa=0$.
Then taking into account that $\psi_0(y)=-\psi_0(-y)$ from (\ref{K1}) we get
$K_1=1$ and $p_1=p_2$ and then from (\ref{K2}) we get $K_2=1$.

Thus, for $p_1\ne p_2$ the Hamiltonians $H_p$ (\ref{s2}) cannot be obtained one from another
by a  diffeomorphismus $\phi :\ S^2\to S^2$.

\hfill $\Box$

\medbreak
\noindent
{\bf Acknowledgement.} 
The author would like to thank Professor Gerhard  Huisken and the Arbeitsbereich 
"Analysis" of the University of  T\"ubingen 
 for their hospitality.
\medbreak
\ \

\bigbreak
\noindent
After the paper was finished, K. P.   Hadeler informed the 
author that the differential equation in (\ref{inipr}) 
can be solved by succsessful 
integration  and the inverse  function of the solution 
$\psi_0(y)$ to the initial value problem  (\ref{inipr})
 can be given as:
$$y=y(\psi_0)=\frac{1}{4}\log \frac{(\psi_0^4+1)^{\frac{1}{4}}+
\psi_0}{(\psi_0^4+1)^{\frac{1}{4}}-\psi_0}-
\frac{1}{2}\arctan\frac{(\psi_0^4+1)^{\frac{1}{4}}}{\psi_0}.$$

\bigbreak

\bigbreak
\small  {\em Mathematisches Institut, Universit\"at T\"ubingen}   
 
\small  {\em Auf der Morgenstelle 10,  72076  T\"ubingen, Germany} 
 
\small  {\em e-mail: \ \ lena@moebius.mathematik.uni-tuebingen.de} 

\medbreak
\small  {\em Department of Geometry, Nizhny Novgorod State 
Pedagogical University}

\small  {\em 603000 Russia, Nizhny Novgorod, ul. Ulyanova 1}

\end{document}